\documentclass[12pt,a4paper]{article}
\def\ba{\begin{array}}
\def\ea{\end{array}}
\def\be{\begin{equation}}
\def\ee{\end{equation}}
\def\lbl{\label}
\def \rf  {(\ref}

\def\x0{\x_0}
\def\x1{\x_1}

\def\real{{\bf R}}
\def\compl{{\bf C}}
\def\cd{{\cal {D}}}
\def\cg{{\cal {G}}}
\def\tcg{{\tilde{\cal G}}}

\newtheorem{theorem}{Theorem}

\begin{document}
\author{L. \v Snobl
 and L. Hlavat\'y
\thanks{{E--mail: snobl@newton.fjfi.cvut.cz}, {hlavaty@br.fjfi.cvut.cz}} 
\\ {\it Faculty of Nuclear Sciences and Physical Engineering,} 
\\{ \it Czech Technical University,} 
\\ {\it B\v rehov\'a 7, 115 19 Prague 1, Czech Republic}}

\title{Classification of 6--dimensional real Drinfeld doubles}
\date{March 27, 2002}

\maketitle
\bibliographystyle{unsrt}

\abstract{Starting from the classification of real Manin triples we look for 
those that are isomorphic as 6--dimensional Drinfeld doubles i.e. Lie algebras with the  ad--invariant form used 
for construction of the Manin triples. We use several invariants of the Lie algebras to distinguish 
the non--isomorphic structures and give explicit form of maps between Manin triples that are 
decompositions of isomorphic Drinfeld doubles.
The result is a complete list of 6--dimensional real Drinfeld doubles. It consists of 22 classes 
of non--isomorphic Drinfeld doubles. 
}
\vskip 1cm
PACS codes: 02.20.Qs, 02.20.Sv, 11.25.Hf
\newpage\section{Introduction}
In recent years, the study of T--duality in string theory has led to 
discovery of Poisson--Lie T--dual sigma models. Klim\v{c}\'{\i}k and \v{S}evera 
have found a procedure allowing to construct the dual models 
from Manin triples  
$({\cal D},{\cal G},{\tilde{\cal G}})$ , 
i.e. a decompositions of a Lie algebra
${\cal D}$ (it must be even--dimensional) 
into two maximally isotropic subalgebras ${\cal G},{\tilde{\cal G}}$ w.r.t. a bilinear form.
The construction of the Poisson--Lie T--dual sigma models is described in \cite{klse:dna} and \cite{kli:pltd}. 

The Lie group possessing a Lie algebra that can be written as a Manin triple is called the Drinfeld double. The classification of the two--dimensional Drinfeld doubles is trivial and the 
four--dimensional Drinfeld doubles can 
be found e.g. in the paper \cite{hlasno:pltdm2dt} together with
the corresponding two-dimensional T--dual 
models. Examples of six--dimensional Drinfeld doubles and three--dimensional dual models  were given e.g. in \cite{vall:su2,sfe:pltd,jare:pltd}. 
There was an attempt to classify the six--dimensional Drinfeld doubles by the Bianchi forms of their 
three--dimensional isotropic subalgebras in \cite{jare:pltd} but it is not 
sufficient for the specification of the Drinfeld double.

As we shall see Manin triples are equivalent to Lie bialgebras and the classification of the three--dimensional Lie bialgebras 
(i.e. six--dimensional Manin triples) was given in \cite{gom:ctd}. Without knowledge of this this work we have performed a classification 
ot the six--dimensional Manin triples \cite{hlasno:mtriples}. The consequent comparison proved that the results are identical even though 
we have started from a different description of the three--dimensional algebras and used a completely different method. It means that 
in \cite{hlasno:mtriples} we have done an independent check of  \cite {gom:ctd} and on the other hand expressed the results 
in a different form, namely as Manin triples. 

The goal of 
this paper is to find which of the Manin triples represent decomposition of the same 
(or more precisely isomorphic) Drinfeld doubles. We use the notation of \cite{hlasno:mtriples} because the less compact sorting of 
the triples into parametrized classes turned out more appropriate for the classification. The result is a complete list of 
the real non--isomorphic six--dimensional Drinfeld doubles. 
Let us note that the Drinfeld double is defined not only by its Lie structure  but also by a bilinear form. There are e.g. two classes 
of Drinfeld doubles for $so(1,3)$ as we shall see. 

In the following sections, we firstly recall the definitions 
of Manin triple, Lie bialgebra and Drinfeld double then briefly explain the 
approach we have used to distinguish the non--isomorphic structures. The main result of the paper is the 
classification theorem in the Section \ref{sec3}. Explicit 
forms of maps between Manin triples that are decompositions of the isomorphic Drinfeld doubles are contained in the proof of the theorem.

\section{Manin triples,  Lie bialgebras, Drinfeld doubles}

The Drinfeld double $D$ is defined as a connected Lie group such that its Lie algebra 
$\cd$ equipped by a symmetric ad--invariant nondegenerate bilinear form 
$\langle\, .,.\rangle $ can be decomposed into a pair of 
subalgebras $\cg$, $\tcg$ maximally isotropic w.r.t. $\langle\, .,.\rangle $ and $\cd$ as 
a vector space is the direct 
sum of $\cg$ and $\tcg$. This ordered triple of algebras $(\cd,\cg$,$\tcg)$ is called Manin triple. 

One can see that the dimensions of the subalgebras are equal and that bases 
$\{X_i\}, \{\tilde X^i\},\ i=1,2,3$ in the subalgebras can be chosen so that
\be \langle X_i,X_j\rangle =0,\  \langle X_i,\tilde X^j\rangle =\langle 
\tilde X^j,X_i\rangle =\delta_i^j,\  \langle \tilde X^i,\tilde X^j\rangle =0.\lbl{brackets}\ee
This canonical form of the bracket is invariant with respect to the transformations 
\be X_i'=X_k A^k_i,\ \tilde X^{'j}=(A^{-1})^j_k \tilde X^k. \lbl{tfnb}\ee
The Manin triples that are related by the transformation \rf{tfnb}) are considered isomorphic.
Due to the ad-invariance of $\langle\, .,.\rangle $ the algebraic structure of $\cd$ is determined 
by the structure of the maximally isotropic subalgebras because in the basis $\{X_i\}, \{\tilde X^i\}$ 
the Lie product is given by 
\[ [X_i,X_j]={f_{ij}}^k X_k,\ [\tilde X^i,\tilde X^j]={\tilde {f^{ij}}_k} \tilde X^k,\]
\be [X_i,\tilde X^j]={f_{ki}}^j \tilde X^k +{\tilde {f^{jk}}_i} X_k. \lbl{liebd}\ee

It is clear that to any Manin triple $({\cal D},{\cal G},{\tilde{\cal G}})$   one can construct the 
dual one by 
interchanging $\cg \leftrightarrow \tcg$, i.e. interchanging the structure coefficients 
$ {f_{ij}}^k \leftrightarrow {\tilde {f^{ij}}_k}$. All properties of Lie algebras 
(the nontrivial being the Jacobi identities) remain to be satisfied. On the other hand for given Drinfeld double more than two 
Manin triples can exist and we shall see many examples of that.
The Drinfeld double $D$ is defined as a Lie group such that its Lie algebra 
$\cd$ equipped by a symmetric ad--invariant nondegenerate bilinear form 
$\langle .,.\rangle $ can be decomposed into a pair of maximally isotropic 
subalgebras $\cg$, $\tcg$ such that $\cd$ as a vector space is the direct 
sum of $\cg$ and $\tcg$. This ordered triple of algebras $(\cd,\cg$,$\tcg)$ is called Manin triple. 

One can rewrite the structure 
of a  Manin triple also in another,
equivalent, but for certain considerations more suitable, form of Lie bialgebra. 

A Lie bialgebra is a Lie algebra $g$ equipped also by a Lie cobracket\footnote{Summation index is suppressed} 
 $\delta:g \rightarrow g \otimes g: \delta(x)=  \sum x_{[1]} \otimes x_{[2]}$ such that 
 \begin{eqnarray}
\sum x_{[1]} \otimes x_{[2]} & = & - \sum x_{[2]} \otimes x_{[1]}, \\
 (id \otimes \delta) \circ \delta (x) & + & {\rm cyclic \, permutations \, of \, tensor \, indices} = 0, \label{dji1} \\
\nonumber \delta([x,y]) & = & \sum [x,y_{[1]}] \otimes y_{[2]} + y_{[1]} \otimes [x,y_{[2]}] - \\
 & - & [y,x_{[1]}] \otimes x_{[2]} - x_{[1]} \otimes [ y, x_{[2]} ] \label{mji1}
\end{eqnarray}
(for detailed account on Lie bialgebras see e.g. \cite{Drinfeld} or \cite{Majid}, Chapter 8).

The correspondence between a Manin triple and a Lie bialgebra can now be formulated in the following way.
Because both subalgebras $\cg$, $\tcg$ of the Manin triple are of the same dimension
and are connected by nondegenerate pairing, it is natural to consider $\tcg$ as a dual $\cg^{*}$ to 
$\cg$ and to use the Lie bracket in $\tcg$ to define the Lie cobracket in $\cg$; 
$\delta(x) $ is given by $\langle \delta(x), \tilde{y} \otimes \tilde{z} \rangle = 
\langle x, [\tilde{y},\tilde{z}] \rangle,$ $\forall \tilde{y},\tilde{z} \in \cg^{*}$,
 i.e. $\delta(X_{i})= \tilde {f^{jk}_{i}}  X_{j} \otimes X_{k} $.
The Jacobi identities in $\tcg$ 
\be\label{jidd}
\tilde {f^{kl}_{m}}  \tilde {f^{ij}_{l}} + \tilde {f^{il}_{m}}  \tilde {f^{jk}_{l}} + \tilde {f^{jl}_{m}}  \tilde {f^{ki}_{l}} =0
\ee 
are then equivalent to the property of cobracket (\ref{dji1}) and 
the $\tcg$--component of the mixed Jacobi identities \footnote{The Jacobi identities $[X_i,[\tilde X^j,\tilde X^k]]+ 
{\rm cyclic} = 0$ lead to both (\ref{mjid})  (terms proportional to $\tilde X^l$) and (\ref{jidd})  (terms
proportional to $X_{l}$).} 
\be \label{mjid}
 \tilde {f^{jk}}_{l} {f_{mi}}^{l} +   \tilde {f^{kl}}_{m} {f_{li}}^{j} + 
  \tilde {f^{jl}}_{i} {f_{lm}}^{k} +   \tilde {f^{jl}}_{m} {f_{il}}^{k} + 
  \tilde {f^{kl}}_{i} {f_{lm}}^{j}   =0
\ee
are equivalent to (\ref{mji1}).  

From now on, we will use the formulation in terms of Manin triples, 
Lie bialgebra formulation of all results can be easily 
derived from it. We also consider 
only algebraic structure, the Drinfeld doubles as the Lie groups can
be obtained in principle by means of exponential map and usual
theorems about relation between Lie groups and Lie algebras 
apply, e.g. there is a one to one correspondence between (finite--dimensional) Lie 
algebras and connected and simply connected Lie groups. The group
structure of the Drinfeld double can be deduced e.g. by taking 
matrix exponential of adjoint representation of its algebra.

We shall consider two Drinfeld doubles isomorphic if they have isomorphic algebraic structure and there is an isomorphism 
transforming one ad-invariant bilinear form to the other. As mentioned above we can always choose a basis  so that the 
bilinear form have canonical form \rf{brackets}) and the Lie product is then given by \rf{liebd}). The Drinfeld doubles 
$\cd$ and $\cd'$ with these special bases $Y_a=(X_1,X_2,X_3,\tilde X^1,\tilde X^2,\tilde X^3)$, $Y'_a=(X'_1,X'_2,X'_3,\tilde X'^1,
\tilde X'^2,\tilde X'^3)$ are isomorphic iff there is an invertible $6\times6$ matrix ${C_a}^b$ such that the linear map given by
\be Y'_a = {C_a}^b Y_b \label{isodd}\ee
transforms the Lie multiplication of $\cd$ into that of $\cd'$ and preserves the canonical form of the bilinear form $\langle.
\, ,.\rangle$. 
This is equivalent to 
\be {C_a}^p {C_b}^q B_{pq}= B_{ab},\  {C_a}^p {C_b}^q {F_{pq}}^r ={ F'_{ab}}^c {C_c}^r \lbl{cpodm}\ee
where ${F_{ab}}^c,\ {F'_{ab}}^c,\ a,b,c=1,\ldots,6$ are structure coefficients of the doubles $\cd$ and $\cd'$ and 
\be B=\left( \begin{array}{cccccc}
0&0&0&1&0&0 \\0&0&0&0&1&0 \\0&0&0&0&0&1 
\\1&0&0&0&0&0 \\0&1&0&0&0&0 \\0&0&1&0&0&0     \end{array}
\right). \lbl{bmat} \ee

\section{Method and result of classification}\lbl{sec3}

As mentioned in the Introduction, there are 78 non--isomorphic classes of Manin triples 
\cite{hlasno:mtriples}.  If we take into account the duality transformation $({\cal D},
{\cal G},{\tilde{\cal G}})\mapsto ({\cal D},{\tilde{\cal G}},{\cal G})$ the number 
is reduced to 44. Their explicit form is given in the Appendix B. 
It follows from \rf{brackets}) and \rf{liebd}) that the structure of the Manin triple can be given 
by the structure coefficients $f_{ij}^k,\, {\tilde {f^{ij}}_k}$ of $\cg$ and $\tcg$ in the special 
basis where relations \rf{brackets}) hold. That's why we usually denote the Manin triples 
$(\cd,\cg,\tcg)$ by $(\cg|\tcg)$ or $(\cg|\tcg|b)$ when a scaling parameter $b$ occurs in the definition 
of the Lie product. Let us note that $(\cg|\tcg|b)$ and $(\cg|\tcg|b')$ are isomorphic up to rescaling of 
$\langle .,.\rangle $.

\begin{table}[t]
\begin{tabular}[c]{|c|c|c|c|l|}
\hline
Signature & Dim. of     & Dim. of           & Dim. of           &   \\ 
of  $K$   & $[\cd,\cd]$ & ${\cd}^2,{\cd}^3$ & ${\cd}_2,{\cd}_3$ &  Manin triples  \\ \hline \hline
(3,3,0)   & 6           & 6,6               & 6,6               & ($9|5|b$), ($8|5.ii|b$), \\ 
          &             &                   &                   & ($7_a|7_{1/a}|b$), ($7_0|5.ii|b$)   \\ \hline
(4,2,0)   & 6           & 6,6               & 6,6               & ($8|5.i|b$),  ($6_a|6_{1/a}.i|b$),  \\ 
          &             &                   &                   & ($6_0|5.iii|b$)  \\ \hline
(0,3,3)   & 6           & 6,6               & 6,6               & ($9|1$) \\ \hline
(2,1,3)   & 6           & 6,6               & 6,6               & ($8|1$), ($8|2.iii$), ($7_0|4|b$),  \\
          &             &                   &                   & ($7_0|5.i$), ($6_0|4.i|b$), ($6_0|5.i$),\\ 
          &             &                   &                   & ($5|2.ii$), ($4|2.iii|b$), \\ \cline{2-5}
          & 3           & 3,3               & 3,3               & ($3|3.i$) \\ \hline
(1,0,5)   & 5           & 5,5               & 1,0               & ($7_{a}|1$), ($7_{a}|2.i$), ($7_{a}|2.ii$), $a>1$\\
          &             &                   &                   & ($6_a|1$), ($6_a|2$), ($6_a|6_{1/a}.ii$), \\ 
          &             &                   &                   & ($6_a|6_{1/a}.iii$), ($6_0|1$), ($6_0|2$), \\
          &             &                   &                   & ($6_0|4.ii$), ($6_0|5.ii$), ($5|1$), ($5|2.i$), \\ 
          &             &                   &                   & ($4|1$), ($4|2.i$), ($4|2.ii$)\\ \cline{2-5}
          & 3           & 3,3               & 1,0               & ($3|1$), ($3|2$), ($3|3.ii$), ($3|3.iii$) \\ \hline 
(0,1,5)   & 5           & 5,5               & 1,0               & ($7_{a}|1$), ($7_{a}|2.i$), ($7_{a}|2.ii$), $a<1$ \\
          &             &                   &                   & ($7_0|1$), ($7_0|2.i$), ($7_0|2.ii$) \\ \hline
(0,0,6)   & 5           & 5,5               & 1,0               & ($7_{a}|1$), ($7_{a}|2.i$), ($7_{a}|2.ii$), $a=1$  \\ \cline{2-5}
          & 3           & 0,0               & 0,0               & ($2|2$) \\ \cline{3-5}
          &             & 2,0               & 0,0               & ($2|2.i$), ($2|2.ii$)  \\ \cline{2-5}
          & 0           & 0,0               & 0,0               & ($1|1$) \\  
\hline
\end{tabular}
\caption{Invariants of Manin triples } \label{tab1}
\end{table}

It is clear that a
direct check which of 44 Manin triples are decomposition of isomorphic Drinfeld 
doubles is a tremendous task. That's why we first evaluate as many invariants of the 
algebras as possible and then sort them into smaller subsets according to the 
values of the invariants. It is clear that only the Manin triples in these subsets can
be decomposition of the same Drinfeld double.
The invariants we have used are: 
\begin{itemize}
\item signature (numbers of positive, negative and zero eigenvalues) of the Killing form, 
\item dimensions of the comutant $[\cd,\cd] \equiv {\cd}^1 \equiv {\cd}_1$ and subalgebras 
created by the repeated Lie multiplication ${\cd}^{i+1}=[{\cd}^i,\cd]$, 
(up to $i=3$, it turns out that for $i \geq 3$ 
${\cd}^{i+1}={\cd}^{i}$). (We have for completeness determined also dimensions of ${\cd}_{i+1}=[{\cd}^i,{\cd}^i]$, 
but they doesn't lead to refinement of our partition.) 

\end{itemize}
The partition of the list of Manin triples according to the values of invariants is
in the Table \ref{tab1}.
The final distinction between non--isomorphic Drinfeld doubles and their decomposition into Manin triples provides the following theorem.

\begin{theorem}  \label{veta}
Any 6--dimensional real Drinfeld double belongs just to one of the following 22 classes and allows 
decomposition into all Manin triples listed in the class and their duals $(\cg \leftrightarrow \tcg)$.
If the class contains parameter $a$ or $b$, the Drinfeld doubles with different values of this parameter 
are non--isomorphic. 
\begin{enumerate}
\item $(9|5|b)$ $\cong $ $(8|5.ii|b)$ $\cong $ $(7_0|5.ii|b)$, $b>0$,
\item $(8|5.i|b)$ $\cong $ $(6_0|5.iii|b)$, $b>0$,
\item $(7_a|7_{1/a}|b)$ $\cong $ $(7_{1/a}|7_a|b)$, $a \geq 1,b \in \real - \{ 0 \},$
\item $(6_a|6_{1/a}.i|b)$ $\cong $ $(6_{1/a}.i|6_a|b)$, $a>1, b \in \real - \{ 0 \},$
\item $(9|1)$,
\item $(8|1)$ $\cong$ $(8|2.iii)$ $\cong$ $(7_0|5.i)$ $\cong$ $(6_0|5.i)$ $\cong$ $(5|2.ii)$,
\item $(7_0|4|b)$  $\cong $ $(4|2.iii|b)$ $\cong$ $(6_0|4.i|-b)$, $b \in \real - \{ 0 \},$  
\item $(3|3.i)$,
\item $(7_{a}|1)$ $\cong $ $(7_{a}|2.i)$ $\cong $ $(7_{a}|2.ii)$, $a>1$,
\item $(6_a|1)$ $\cong $ $(6_a|2)$ $\cong $ $(6_a|6_{1/a}.ii)$ $\cong $ $(6_a|6_{1/a}.iii)$, $a>1,$
\item $(6_0|1)$ $\cong$ $(6_0|5.ii)$ $\cong$ $(5|1)$ $\cong$ $(5|2.i)$,
\item $(6_0|2)$ $\cong$ $(6_0|4.ii)$ $\cong$ $(4|1)$ $\cong$ $(4|2.i)$ $\cong$ $(4|2.ii)$,
\item $(3|1)$ $\cong$ $(3|2)$ $\cong$ $(3|3.ii)$ $\cong$ $(3|3.iii)$,
\item $(7_{a}|1)$ $\cong $ $(7_{a}|2.i)$ $\cong $ $(7_{a}|2.ii)$, $0<a<1$,
\item $(7_0|1)$, 
\item $(7_0|2.i)$,
\item $(7_0|2.ii)$,
\item $(7_{1}|1)$ $\cong $ $(7_{1}|2.i)$ $\cong $ $(7_{1}|2.ii)$,
\item $(2|1)$,
\item $(2|2.i)$,
\item $(2|2.ii)$,
\item $(1|1)$.
\end{enumerate}
\end{theorem}

\section{The proof of Theorem \ref{veta}}
The essence of the proof is to find which of the 78 non--isomorphic Manin triples found in \cite{hlasno:mtriples} and displayed in 
the Appendix B yield isomorphic Drinfeld doubles. The isomorphisms are  given by the explicit form of the transformation matrices 
$C$ (see \rf{isodd}) ) that were found by solution of equations  \rf{cpodm}). In this part we have used the computer programs 
Maple V and Mathematica 4. The solutions are not unique and here we present only a simple examples of them. The non--isomorphic
 Drinfeld doubles are distinguished by investigation of their various subalgebras and properties of $\langle\, .,.\rangle $ and 
the Killing form on them.

In the next subsection we analyze the subsets of non--isomorphic Manin triples characterized by invariants described in the Section 
\ref{sec3} and displayed in the Table \ref{tab1}.

\subsection{Manin triples with the Killing form of signature (3,3,0)}
In this case the signature of the Killing form itself fixes the Lie algebra $\cd$ 
of the Drinfeld double uniquely.
It is the well-known $so(3,1)$ which is simple as a real Lie algebra and its 
complexification is semisimple; it decomposes into two copies of $sl(2,{\compl})$. 
The Drinfeld doubles corresponding  to ($9|5|b$), ($8|5.ii|b$), ($7_0|5.ii|b$), ($7_a|7_{1/a}|b$)
can consequently differ only by the bilinear form $\langle \,.,. \rangle$.

We can find a necessary condition for equivalence of semisimple Drinfeld doubles 
from the fact that any invariant symmetric bilinear form on a complex 
simple Lie algebra is a multiple of the Killing form and that any invariant 
symmetric bilinear form on a semisimple Lie algebra is a sum of 
invariant symmetric bilinear forms on its simple components.  
({\it Proof:} Let $\cg = \oplus_{i} \cg_i$ be the decomposition into simple components,
$X \in \cg_i, Y \in \cg_j, \ i \neq j$. Then $\exists A_k, B_k \in \cg_j$ s.t. $Y=\sum_k [A_k, B_k]$ 
and from the ad--invariance of the form $\langle X, Y \rangle = \sum_k \langle X, [A_k, B_k] \rangle =
- \sum_k \langle [A_k,X], B_k \rangle = - \sum_k \langle 0, B_k \rangle =0$.)

We therefore consider the complexification $\cd_{\compl}$ of the Drinfeld double algebra 
and write both the Killing form on $\cd_{\compl}$ and the bilinear form
$\langle \,.,. \rangle$ in terms of Killing forms $K_1,K_2$ 
of still unspecified simple components $sl(2,{\compl})_1$, $sl(2,{\compl})_2$ 
($\cd_{\compl} = sl(2,{\compl})_1 \oplus sl(2,{\compl})_2$)
$$ K= K_1+K_2, \; \langle , \rangle =\alpha K_1 + \beta K_2.$$
We trivially extend the Killing forms $K_1,K_2$ to the whole Drinfeld double 
algebra $\cd_{\compl}$ and express them  as
$$ K_1= \frac{\langle , \rangle-\beta K}{\alpha-\beta}, \, K_2 = \frac{\alpha K - 
\langle , \rangle}{\alpha-\beta}.$$
Because $K_1,K_2$ are trivially extended Killing forms, they  must
have 3--dimensional nullspace ($sl(2,{\compl})_2$ in the case of $K_1$ and 
$sl(2,{\compl})_1$ in the case of $K_2$). These two conditions on dimensions 
of nullspaces fix the coefficients $\alpha, \beta$
uniquely up to a permutation. Therefore, {\bf 
the necessary condition for equivalence of two semisimple 
6--dimensional Drinfeld doubles 
is the equality of their sets of coefficients $\{ \alpha, \beta \} $ .}

We compute  the coefficients $\alpha, \beta$ for the Manin triples in this class
 and find that in three cases ($9|5|b$), ($8|5.ii|b$), ($7_0|5.ii|b$) is
$$\{ \alpha, \beta \} = \{ \frac{i}{4b},-\frac{i}{4b} \} $$
and for ($7_a|7_{1/a}|b$) is $$ \{ \alpha, \beta \} = \{ \frac{ia}{4b(a-i)^2},-\frac{ia}{4b(i+a)^2} \}. $$

We see that the Manin triple ($7_a|7_{1/a}|b$) defines for any $a,b$  Drinfeld doubles different from any of the Drinfeld doubles
associated to the Manin triples ($9|5|b$), ($8|5.ii|b$), ($7_0|5.ii|b$) and that Drinfeld doubles corresponding to ($7_a|7_{1/a}|b$) 
with different values of $a$ and $b$ are different except the case $a'=1/a, b'=b$. The  Manin triples 
($7_a|7_{1/a}|b$) and ($7_{1/a}|7_a|b$) are mutually dual, correspond to  $\cg \leftrightarrow \tcg$ and therefore give the same 
Drinfeld double. The Manin triple ($7_1|7_{1}|b$) is of course self--dual.

Also one sees that the Manin triples 
($9|5|b$), ($8|5.ii|b$), ($7_0|5.ii|b$) with 
different $b$ cannot lead to the same Drinfeld double. For the Manin triples ($9|5|b$), ($8|5.ii|b$), ($7_0|5.ii|b$) with equal $b$, 
the transformations (\ref{isodd}) between Drinfeld doubles exist, but may contain complex numbers since 
 up to now we have considered only complexifications of the original Manin triples.

However, one can check that the following real transformation matrices $C$ guarantee
the equivalence of the Drinfeld doubles in this class for fixed value of $b$.
\[ (9|5|b)\rightarrow (8|5.ii|b)\ : \ \ C= 
\left( \begin{array}{rrrrrr} 
0& 1& 0& 0& 0& \frac{1}{b} \\
 0& 0& 1& 0& -\frac{1}{b}& 0 \\
 1& 0& 0& 0& 0& 0 \\ 
0& 0& 0& 0& 1& 0 \\
 0& 0& 0& 0& 0& 1 \\
 0& 0& 0& 1& 0& 0
\end{array} \right),
  \]
\[ (9|5|b)\rightarrow (7_0|5.ii|b)\ : \ \ C= 
\left( \begin{array}{rrrrrr} 
\frac{1}{2}& 0& -\frac{1}{2}& 0& \frac{1}{2b}& 0 \\
 0& \frac{1}{2}& 0& -\frac{1}{2b}& 0& 0 \\
 0& 0& 1& 0& 0& 0 \\
 0& b& 0& 1& 0& 0 \\
 -b& 0& -b& 0& 1& 0 \\
 0& b& 0& 0& 0& 1
\end{array} \right).
\]
As mentioned in the beginning of this section the transformation matrices are not unique; they contain several free parameters. 
Here and further we give them in a simple form setting the parameters zero or one.

\subsection{Manin triples with the Killing form of signature (4,2,0)}
In this case the signature of the Killing form again fixes the Lie algebra $\cd$ 
of the Drinfeld double uniquely,
it is $sl(2,\real) \oplus sl(2,\real)$, and the Drinfeld doubles  
may again differ only by the bilinear form $\langle \,.,. \rangle$. We use the criterion developed in 
the previous subsection for semisimple Drinfeld doubles and find
\begin{itemize}
\item ($8|5.i|b$), ($6_0|5.iii|b$): $ \{ \alpha, \beta \} = \{ \frac{1}{4b}, -\frac{1}{4b} \} $,
\item ($6_a|6_{1/a}.i|b$):  $ \{ \alpha, \beta \} = \{ \frac{a}{4b(a-1)^2}, - \frac{a}{4b(1+a)^2} \} $.
\end{itemize}
This shows that the Manin triples might specify isomorphic Drinfeld doubles only in the following two cases:
\begin{enumerate}
\item ($8|5.i|b$) and  ($6_0|5.iii|b$) for the same value of $b$. In this case we have found 
the transformation matrix $C$
\[ (8|5.i|b) \rightarrow (6_0|5.iii|b) \  : \ \ C= 
\left( \begin{array}{rrrrrr} 
0 &  0 &  -\frac{b}{2} &  -\frac{1}{2} &  0 &  0 \\
 -\frac{b}{2} &  \frac{b}{2} &  0 &  0 &  0 &  \frac{1}{2} \\ 
0 &  -1 &  0 &  0 &  0 &  0 \\
 -1 &  -1 &  0 &  0 &  0 &  -\frac{1}{b} \\
 0 &  0 &  1 &  -\frac{1}{b} &  0 &  0 \\
 0 &  0 &  b &  0 &  -1 &  0
\end{array} \right).
\]
This transformation is real and therefore the Drinfeld doubles are isomorphic, 
($8|5.i|b$) $\cong $ ($6_0|5.iii|b$). 
\item ($6_a|6_{1/a}.i|b$) and ($6_{1/a}|6_{a}.i|b$). 
One can easily see that these Manin 
triples are dual (i.e. can be obtained one from the other by the interchange $\cg \leftrightarrow \tcg$) 
and the Drinfeld doubles are therefore isomorphic.
\end{enumerate}

\subsection{Manin triples with the Killing form of signature (0,3,3)}
This class contains only one Manin triple ($9|1$) and its dual; the corresponding Drinfeld double is isomorphic 
to $so(3) \triangleright  {\real}^3  $ since the Killing form has the signature (0,3,3) and ${\rm dim} [\cd,\cd]=3$.

\subsection{Manin triples with the Killing form of signature (2,1,3)}
We consider only the Manin triples with ${\rm dim}[\cd, \cd]=6$, the other set in this class 
contains only one Manin triple ($3|3.i$), which is isomorphic as a Lie algebra to $sl(2,\real) \oplus {\real}^3$ 
since the Killing form has the signature (2,1,3) and ${\rm dim} [\cd,\cd]=3$.

The Manin triples in this set 
($8|1$), ($8|2.iii$), ($7_0|4|b$), ($7_0|5.i$), ($6_0|4.i|b$), ($6_0|5.i$),  
($5|2.ii$), ($4|2.iii|b$), are neither semisimple (${\rm rank} K \neq 6$) nor solvable 
($[\cd,\cd]= \cd$). Therefore they have a nontrivial Levi--Maltsev decomposition into 
semidirect sum of a semisimple subalgebra $S$ and radical $N$
$$ \cd = S \triangleright N,$$ 
both of them are 3--dimensional.
 Knowledge of this decomposition turns out to be helpful in the investigation 
of equivalence of the Drinfeld doubles.

A rather simple computation shows that the radical  is in all these Manin triples abelian 
and maximally isotropic, e.g. for ($8|1$) the radical is $N={\rm span} \{ \tilde{X}^1, 
\tilde{X}^2, \tilde{X}^3 \} $, for ($4|2.iii|b$) the radical is $N={\rm span} 
\{ X_3 , \tilde{X}^1, \tilde{X}^2 \} $.

Next we find the semisimple component. It turns out that the semisimple subalgebra $S$
is in all cases $sl(2,\real)$, e.g. for ($8|1$) it can be evidently chosen $S = {\rm span} \{ {X}_1, 
{X}_2, {X}_3 \} $, for ($4|2.iii|b$) the most general form of the semisimple subalgebra is 
$S={\rm span} \{  
2 X_1 -2\alpha X_3 -\frac{2}{b} \tilde{X}^1 -2 \beta \tilde{X}^2,
-\frac{2}{b} X_2 -\frac{2\gamma}{b} X_3 - \frac{2 \beta}{b} \tilde{X}^1,
 \alpha \tilde{X}^1 + (2-\gamma) \tilde{X}^2 + \tilde{X}^3
\} $ for any values of $\alpha,\beta,\gamma$.

One can restrict the form $\langle\, . , .\rangle$ to the semisimple subalgebra $S$ and finds that for 
 ($8|1$) $\langle\, . , .\rangle_S=0$, i.e. $S$ is maximally isotropic, whereas for ($4|2.iii|b$) 
and any choice of $\alpha,\beta,\gamma$ is 
$\langle\, . , .\rangle_S = - 1/ b K_S$, $K_S$ being the Killing form on $S$. This shows that 
as Drinfeld doubles ($8|1$) and ($4|2.iii|b$) and similarly ($4|2.iii|b$) for different values of $b$ 
are not isomorphic.

Performing the same computation for all Manin triples in this set, we find that they divide into two
subsets. 
\begin{enumerate}
\item ($8|1$), ($8|2.iii$),  ($7_0|5.i$),  ($6_0|5.i$),  ($5|2.ii$): $\langle\, . , .\rangle_S=0$
\item ($7_0|4|b$), ($6_0|4.i|-b$), ($4|2.iii|b$): $\langle\, . , .\rangle_S = - 1/b K_S$, $b \in \real - \{
0 \} $
\end{enumerate}
We find the transformation matrices for Manin triples in these subsets and prove the equivalence of
the corresponding Drinfeld doubles:
\[ (8|1) \rightarrow (8|2.iii) \ : \ \ C= 
\left( \begin{array}{rrrrrr} 
 -1& 0& 0& 0& 0& 0 \\ 0& -1& 0& 0& 0& 0 \\ 0& 0& 1& 0& 0& 0 \\ 
0& 1& 0& -1& 0& 0 \\ -1& 0& -1& 0& -1& 0 \\ 0& -1& 0& 0& 0& 1 
\end{array} \right),
  \]
\[ (8|1) \rightarrow  (7_0|5.i) \ : \ \ C= 
\left( \begin{array}{rrrrrr} 
 0& 0& 0& 0& -1& 0 \\ 0& 0& 0& 1& 0& 0 \\ 0& 0& 1& 0& 0& 0 \\ 
0& -1& 0& 0& 0& 0 \\ 1& 0& 1& 0& 0& 0 \\ 0& 0& 0& -1& 0& 1 
\end{array} \right),
  \]
\[ (8|1) \rightarrow  (6_0|5.i) \ : \ \ C= 
\left( \begin{array}{rrrrrr} 
 0& 0& 0& 0& 1& 0 \\ 0& 0& 0& 0& 0& -1 \\ 1& 0& 0& 0& 0& 0 \\ 
0& 1& 0& 0& 0& 0 \\ 1& 0& -1& 0& 0& 0 \\ 0& 0& 0& 1& 0& 1 
\end{array} \right),
  \]
\[ (8|1) \rightarrow   (5|2.ii) \ : \ \ C= 
\left( \begin{array}{rrrrrr} 
 1& -1& -1& -1& -1& 0 \\ 0& 0& 0& 0& -1& 1 \\ 0& -1& -1& -1& 0& 0 \\ 
0& 0& 0& 1& -1& 1 \\ -1& 0& 1& 0& 1& 0 \\ 0& 0& 0& -1& 0& -1 
\end{array} \right),
  \]
respectively
\[  (4|2.iii|b) \rightarrow (7_0|4|b) \ : \ \ C= 
\left( \begin{array}{rrrrrr} 
 0& 0& 0& \frac{1}{b} & 0& 0 \\ 0& 0& -\frac{1}{2b} & 0& 1& 0 \\ 
0& \frac{1}{2b} & 0& 0& 0& 1 \\ b & 0& 0& 0& 0& 0 \\ 0& 1& 0& 0& 0& 0 \\ 0& 0& 1& 0& 0& 0 
\end{array} \right),
  \]
\[  (4|2.iii|b) \rightarrow (6_0|4.i|-b) \ : \ \ C= 
\left( \begin{array}{rrrrrr} 
 0& 0& 0& -\frac{1}{b} & 0& 0 \\ 0& 0& \frac{1}{2b} & 0& 1& 0 \\ 
0& -\frac{1}{2b} & 0& 0& 0& 1 \\ -b & 0& 0& 0& 0& 0 \\ 0& 1& 0& 0& 0& 0 \\ 0& 0& 1& 0& 0& 0 
\end{array} \right).
  \]

Concerning the Lie structure of these Drinfeld doubles, it follows from the signature of the Killing form 
and dimension of $[\cd,\cd]$ that the Lie algebra of $D$ is isomorphic in both cases to $sl(2,{\real}) 
\triangleright{\real}^3$ where commutation relations between subalgebras are given by the unique
irreducible representation of $sl(2,{\real})$ on ${\real}^3$.

\subsection{Manin triples with the Killing form of signature (1,0,5)}
\subsubsection{Case ${ \rm \bf dim}{\bf [\cd, \cd]=5}$}\label{nilp105}
This set contains the greatest number of Manin triples:
($7_{a>1}|1$),  ($7_{a>1}|2.i$), ($7_{a>1}|2.ii$), 
($6_a|1$),  ($6_a|2$),  ($6_a|6_{1/a}.ii$),  ($6_a|6_{1/a}.iii$),
($6_0|1$),  ($6_0|5.ii$),  ($5|1$),  ($5|2.i$),
($6_0|2$),  ($6_0|4.ii$), ($4|1$),  ($4|2.i$),  ($4|2.ii$).
In order to shorten our considerations we firstly present the transformation matrices $C$ 
showing the equivalence of following Drinfeld doubles and later we prove that the following classes 
of Drinfeld doubles are non--isomorphic:
\begin{enumerate}
\item ($7_{a>1}|1$) $\cong $ ($7_{a>1}|2.i$) $\cong $ ($7_{a>1}|2.ii$) for the same value of $a$ 
\[  (7_{a}|1) \rightarrow (7_{a}|2.i) \ : \ \ C= 
\left( \begin{array}{rrrrrr} 
 1 &  0 &  0 &  0 &  0 &  0 \\
 0 &  1 &  0 &  0 &  0 &  0 \\
 0 &  0 &  1 &  0 &  0 &  0 \\ 
0 &  0 &  0 &  1 &  0 &  0 \\
 0 &  0 &  -\frac{1}{2a} &  0 &  1 &  0 \\
 0 &  \frac{1}{2a} &  0 &  0 &  0 &  1  
\end{array} \right) ,  
  \]
\[  (7_{a}|1) \rightarrow (7_{a}|2.ii) \ : \ \ C= 
\left( \begin{array}{rrrrrr} 
 -1 &  0 &  0 &  0 &  0 &  0 \\
 0 &  0 &  0 &  0 &  -2a &  0 \\
 0 &  0 &  0 &  0 &  0 &  2a \\ 
0 &  0 &  0 &  -1 &  0 &  0 \\
 0 &  -\frac{1}{2a} &  0 &  0 &  0 &  1 \\
 0 &  0 &  \frac{1}{2a} &  0 &  1 &  0  
\end{array} \right) .  
  \]
\item ($6_a|1$) $\cong $ ($6_a|2$) $\cong $ ($6_a|6_{1/a}.ii$) $\cong $ ($6_a|6_{1/a}.iii$)
 for the same value of $a$ 
\[  (6_a|1) \rightarrow (6_a|2) \ : \ \ C= 
\left( \begin{array}{rrrrrr} 
 1 &  0 &  0 &  0 &  0 &  0 \\
 0 &  1 &  0 &  0 &  0 &  0 \\
 0 &  0 &  1 &  0 &  0 &  0 \\ 
0 &  0 &  0 &  1 &  0 &  0 \\
 0 &  0 &  -\frac{1}{2a} &  0 &  1 &  0 \\
 0 &  \frac{1}{2a} &  0 &  0 &  0 &  1  
\end{array} \right) ,  
  \]
\[  (6_a|1) \rightarrow  (6_a|6_{\frac{1}{a}}.ii) \ : \ \ C= 
\left( \begin{array}{rrrrrr} 
 1 &  0 &  0 &  0 &  0 &  1 \\
 0 &  0 &  1-a &  a-1 &  0 &  0 \\
 0 &  1-a &  0 &  0 &  0 &  0 \\ 
0 &  -1 &  1 &  0 &  0 &  0 \\
 \frac{1}{a-1} &  0 &  0 &  0 &  0 &  0 \\
 -\frac{1}{a-1} &  0 &  0 &  0 &  -\frac{1}{a-1} &  -\frac{1}{a-1}  
\end{array} \right) ,  
  \]
\[  (6_a|1) \rightarrow  (6_a|6_{\frac{1}{a}}.iii) \ : \ \ C= 
\left( \begin{array}{rrrrrr} 
 1 &  0 &  0 &  0 &  0 &  1 \\
 0 &  0 &  -1-a &  a+1 &  0 &  0 \\
 0 &  -1-a &  0 &  0 &  0 &  0 \\ 
0 &  1 &  1 &  0 &  0 &  0 \\
 \frac{1}{a+1} &  0 &  0 &  0 &  0 &  0 \\
 \frac{1}{a+1} &  0 &  0 &  0 &  -\frac{1}{a+1} &  \frac{1}{a+1}  
\end{array} \right) .  
  \]
\item ($5|1$) $\cong$ ($5|2.i$) $\cong$  ($6_0|1$) $\cong$ ($6_0|5.ii$) 
\[  (5|1) \rightarrow  (5|2.i) \ : \ \ C= 
\left( \begin{array}{rrrrrr} 
  -1 & 0 & 0 & 0 & 0 & 0 \\
  0 & 0 & 0 & 0 & 1 & 0 \\
  0 & 0 & 0 & 0 & 0 & 1 \\
  0 & 0 & 0 & -1 & 0 & 0 \\
  0 & 1 & 0 & 0 & 0 & - \frac{1}{2}   \\
  0 & 0 & 1 & 0 & \frac{1}{2} & 0  
\end{array} \right) ,  
  \]
\[  (5|1)    \rightarrow (6_0|1)   \ : \ \ C= 
\left( \begin{array}{rrrrrr}
  0 & 0 & - \frac{1}{2}   & 0 & 1 & 0 \\
  0 & 0 & \frac{1}{2} & 0 & 1 & 0 \\
  -1 & 0 & 0 & 0 & 0 & 0 \\
  0 & \frac{1}{2} & 0 & 0 & 0 & -1 \\
  0 & \frac{1}{2} & 0 & 0 & 0 & 1 \\
  0 & 0 & 0 & -1 & 0 & 0 
\end{array} \right) ,  
  \]
\[  (5|1)    \rightarrow (6_0|5.ii)   \ : \ \ C= 
\left( \begin{array}{rrrrrr} 
  0 & -1 & 0 & 0 & 0 & \frac{1}{2} \\
  0 & 1 & 0 & 1 & 0 & \frac{1}{2} \\
  -1 & 0 & 1 & 0 & \frac{1}{2} & 0 \\
  1 & 0 & 0 & 0 & -1 & 0 \\
  1 & 0 & 0 & 0 & 0 & 0 \\
  0 & 0 & 0 & 0 & 0 & 1  
\end{array} \right) .  
  \]
\item ($4|1$) $\cong$ ($4|2.i$) $\cong$ ($4|2.ii$) $\cong$ ($6_0|2$) $\cong$ ($6_0|4.ii$) 
\[  (4|1) \rightarrow  (4|2.i) \ : \ \ C= 
\left( \begin{array}{rrrrrr} 
1 &0 &0 &0 &0 &0 \\
  0 &1 &0 &0 &0 &0 \\
  0 &0 &1 &0 &0 &0 \\ 
 0 &0 &0 &1 &0 &0 \\
  0 &0 & - \frac{1}{2}   &0 &1 &0 \\ 
 0 &\frac{1}{2} &0 &0 &0 &1   
\end{array} \right) ,  
  \]
\[  (4|1) \rightarrow  (4|2.ii) \ : \ \ C= 
\left( \begin{array}{rrrrrr} 
1 & 0 & 0 & 0 & 0 & 0 \\
 0 & 1 & 0 & 0 & 0 & 0 \\
 0 & 0 & 1 & 0 & 0 & 0 \\ 
 0 & 0 & 0 & 1 & 0 & 0 \\
 0 & 0 & \frac{1}{2} & 0 & 1 & 0 \\ 
 0 & - \frac{1}{2}   & 0 & 0 & 0 & 1
\end{array} \right) ,  
  \]
\[  (4|1) \rightarrow  (6_0|2) \ : \ \ C= 
\left( \begin{array}{rrrrrr} 
 0 &  0 &  \frac{1}{2} &  0 &  1 &  0 \\
 0 &  0 &  \frac{1}{2} &  0 &  -1 &  0 \\
 1 &  0 &  0 &  0 &  0 &  0 \\ 
0 &  \frac{1}{2} &  0 &  0 &  0 &  1 \\
 0 &  -\frac{1}{2} &  0 &  0 &  0 &  1 \\
 0 &  0 &  0 &  1 &  0 &  0  
\end{array} \right)  ,  
  \]
\[  (4|1) \rightarrow  (6_0|4.ii) \ : \ \ C= 
\left( \begin{array}{rrrrrr} 
 0 &  0 &  1 &  1 &  \frac{1}{2} &  0 \\
 0 &  0 &  -1 &  0 &  \frac{1}{2} &  0 \\
 -1 &  1 &  0 &  0 &  0 &  \frac{1}{2} \\ 
1 &  0 &  0 &  0 &  0 &  0 \\
 1 &  0 &  0 &  0 &  0 &  -1 \\
 0 &  0 &  0 &  0 &  1 &  0  
\end{array} \right) .  
  \]
\end{enumerate}

In the proof of inequivalence of the above given classes of Manin triples
  we exploit the fact that the Drinfeld doubles 
have at least one decomposition into Manin triple with the 2nd subalgebra  $\tcg$
abelian; we will use only these representants ($7_{a}|1$),$a>1$,
($6_a|1$), ($5|1$),  ($4|1$) in our considerations.

Firstly we find all maximal isotropic abelian subalgebras ${\cal A}$ 
of each of the given Drinfeld doubles. The dimension of any such ${\cal A}$ must be $3$
from the maximal isotropy.
The commutant is in all these cases 
$\cd_1=[\cd,\cd] = {\rm span} \{ X_2,X_3,\tilde{X}^1,\tilde{X}^2,\tilde{X}^3 \}$ 
and the centre is $Z(\cd)={\rm span}\{ \tilde{X}^1 \}=\cd_2$. 
One can see that any element of the form $X_1+Y, \, Y \in \cd_1$ cannot occur
in  ${\cal A}$ because $X_1$ commutes only with $Z(\cd)$ and itself.
Therefore, ${\cal A} \subset \cd_1$. 
Further it follows from the maximality that  ${\cal A}$ 
contains $Z(\cd)$ and we conclude that  ${\cal A}={\rm span} \{ \tilde{X}^1, Y_1, Y_2 \} $ where
$Y_1, Y_2 \in {\rm span} \{ X_2,X_3,\tilde{X}^2, \tilde{X}^3 \}$. 
Analyzing the maximal isotropy and replacing $Y_1,Y_2$ by their suitable linear combinations  
we find that ${\cal A}$  can be in general expressed in one of the following forms
\begin{enumerate}
\item ${\cal A} = {\rm span} \{ \tilde{X}^1, {X}_2, \tilde{X}^3  \},$
\item ${\cal A} = {\rm span} \{ \tilde{X}^1,X_2+ \alpha \tilde{X}^3, X_3 - \alpha \tilde{X}^2  \},$
\item ${\cal A} = {\rm span} \{ \tilde{X}^1, X_2+\alpha X_3, -\alpha \tilde{X}^2+\tilde{X}^3  \},$
\item ${\cal A} = {\rm span} \{ \tilde{X}^1, {X}_3, \tilde{X}^2 \},$
\item ${\cal A} = {\rm span} \{ \tilde{X}^1, \tilde{X}^2, \tilde{X}^3 \}.$
\end{enumerate}

In the next step we check which of these subspaces really form a subalgebra of the given Manin triple.
\begin{itemize}
\item ($7_{a}|1$): the maximal isotropic abelian subalgebras are 
${\rm span} \{ \tilde{X}^1,X_2, X_3 \}$ and ${\rm span} \{ \tilde{X}^1, \tilde{X}^2, \tilde{X}^3 \}$.
One may easily construct for each of these maximal isotropic abelian subalgebras the dual 
(w.r.t $\langle \, . \rangle$) subalgebra 
by taking the remaining elements of the standard basis $X_1, \ldots, \tilde{X}^3$
and finds that it is isomorphic in both cases to Bianchi algebra $7_{a}$. In other words, 
we have shown that 
this class of Drinfeld doubles is non--isomorphic to the other ones and are mutually non--isomorphic for
different values of $a$.
\item  ($6_a|1$): the maximal isotropic abelian subalgebras are 
${\rm span} \{ \tilde{X}^1,X_2, X_3 \}$,
${\rm span} \{ \tilde{X}^1, X_2+ X_3, - \tilde{X}^2+\tilde{X}^3  \},$
${\rm span} \{ \tilde{X}^1, X_2- X_3, \tilde{X}^2+\tilde{X}^3  \},$
${\rm span} \{ \tilde{X}^1, \tilde{X}^2, \tilde{X}^3 \}$.
By a slightly more complicated construction of the dual subalgebras 
we find that they are of the Bianchi type $6_a$ for the same $a$, i.e. this class 
of Drinfeld doubles is non--isomorphic to the other ones and are mutually non--isomorphic for
different values of $a$.
\item ($5|1$): the maximal isotropic abelian subalgebras are 
${\rm span} \{ \tilde{X}^1, {X}_2, \tilde{X}^3  \},$
${\rm span} \{ \tilde{X}^1, X_2, X_3   \},$
${\rm span} \{ \tilde{X}^1, X_2+\alpha X_3, -\alpha \tilde{X}^2+\tilde{X}^3  \},$
${\rm span} \{ \tilde{X}^1, {X}_3, \tilde{X}^2 \}$ and
${\rm span} \{ \tilde{X}^1, \tilde{X}^2, \tilde{X}^3 \}.$
\item ($4|1$): the maximal isotropic abelian subalgebras are 
$ {\rm span} \{ \tilde{X}^1,X_2+ \alpha \tilde{X}^3, X_3 - \alpha \tilde{X}^2  \},$
${\rm span} \{ \tilde{X}^1, {X}_3, \tilde{X}^2 \}$ and
${\rm span} \{ \tilde{X}^1, \tilde{X}^2, \tilde{X}^3 \}.$
\end{itemize}
Already from comparison of number of possible maximal isotropic abelian subalgebras 
for ($5|1$) and ($4|1$)
one sees that the corresponding Drinfeld doubles are non--isomorphic.

It also follows that Drinfeld doubles corresponding to Manin triples ($7_{a}|1$), ($6_a|1$), ($5|1$) and ($4|1$)
are different as Lie algebras, since any maximal isotropic abelian subalgebra $\cal A$ of these Manin triples 
is in fact an abelian ideal ${\cal I}$ such that $[{\cd},{\cal I}]={\cal I}$ and any such 3--dimensional ideal is 
maximal isotropic from ad--invariance of $\langle \ , \ \rangle$. Therefore we have in fact identified the non--isomorphic Drinfeld doubles 
from the knowledge of these ideals $\cal I$ (and in some cases 
$\cd/{\cal I}$) which doesn't depend on the form $\langle \ , \ \rangle$ and the doubles differ already in their Lie algebra structure.

\subsubsection{Case ${ \rm \bf dim}{\bf [\cd, \cd]=3}$}
All Manin triples of this subset are decomposition of one Drinfeld double, 
i.e. they can be transformed one into another by the transformation \rf{isodd}). 
Below are the corresponding matrices.
\[ (3|1)\rightarrow (3|2)\ : \ \ C= 
\left( \begin{array}{rrrrrr} 
-1 & 0 & 0 & 0 & 0 & 0 
\cr 0 & \frac{1}{2} & - \frac{1}{2} & 0 & 1 & 1 
\cr 0 & - \frac{1}{2}   & \frac{1}{2} & 0 & 1 & 1 
\cr 0 & 0 & 0 & -1 & 0 & 0 \cr
 0 & \frac{1}{2} & 0 & 0 & 0 & -1 \cr
 0 & \frac{1}{2} & 0 & 0 & 0 & 1 
\end{array} \right),
  \]
\[ 
(3|1)\rightarrow (3|3.ii)\ : \ \ C= 
\left( \begin{array}{rrrrrr}  
  1 &0 &0 &0 &0 &2 \\
  0 &1 &0 &0 &0 &0 \\
  0 &0 &1 &-2 &0 &0 \\ 
  0 &- \frac{1}{2}   &\frac{1}{2} &0 &0 &0 \\ 
  \frac{1}{2} &0 &0 &0 &1 &1 \\ 
  -\frac{1}{2}   &0 &0 &0 &0 &0
\end{array} \right) ,
\]
\[ 
(3|1)\rightarrow (3|3.iii)\ : \ \ C= 
\left( \begin{array}{rrrrrr}  
 -1 &0 &0 &0 &0 &0 \\
 0 &0 &0 &0 &0 &-2 \\ 
 0 &0 &0 &0 &-2 &0\\
 0 &0 &0 &-1 &1 &1\\ 
- \frac{1}{2}   &0 & - \frac{1}{2}   &0 &0 &0\\
 - \frac{1}{2}   &- \frac{1}{2}   &0 &0 &0 &0
\end{array} \right).
\]

\subsection{Manin triples with the  Killing form of signature (0,1,5)}
This set contains Manin triples ($7_{a<1}|1$), ($7_{a<1}|2.i$), ($7_{a<1}|2.ii$),
($7_0|1$), ($7_0|2.i$), ($7_0|2.ii$). As in the subsection \ref{nilp105}
we can show that Manin triples ($7_{a<1}|1$), ($7_{a<1}|2.i$), ($7_{a<1}|2.ii$) are
decomposition of  isomorphic Drinfeld doubles for the same $a$; the transformation matrices
given above for $a>1$
are meaningful also in this case. It remains to be investigated whether the  
Drinfeld doubles induced by ($7_0|1$), ($7_0|2.i$), ($7_0|2.ii$) are isomorphic as or not.

We again find all maximal isotropic abelian subalgebras of these Manin triples. 
We find
\begin{itemize}
\item ($7_0|1$): the maximal isotropic abelian subalgebras are \\
${\rm span} \{ \tilde{X}^3, X_1+\alpha \tilde{X}^2,  X_2 - \alpha \tilde{X}^1  \},$
$ {\rm span} \{ \tilde{X}^1, \tilde{X}^2, \tilde{X}^3   \},$
\item ($7_0|2.i$) the only maximal isotropic abelian subalgebra is \\
${\rm span} \{ \tilde{X}^3, X_1,  X_2   \},$ the dual subalgebra to it w.r.t 
$\langle \, .,.\rangle$  doesn't exist. 
\item ($7_0|2.ii$) the only maximal isotropic abelian subalgebra is \\
${\rm span} \{ \tilde{X}^3, X_1,  X_2   \},$ the dual subalgebra to it w.r.t 
$\langle \, .,.\rangle$  doesn't exist.
\end{itemize}
This means that Drinfeld double induced by ($7_0|1$) has only decompositions into Manin triple ($7_0|1$) 
and that Drinfeld doubles corresponding to ($7_0|2.i$), ($7_0|2.ii$) are not isomorphic to 
the Drinfeld double corresponding to  ($7_{a<1}|1$)
for any value of $a$. To prove that also ($7_0|2.i$), ($7_0|2.ii$) induce non--isomorphic Drinfeld 
doubles, we find all isotropic subalgebras of Bianchi type $7_0$ in the Manin triple ($7_0|2.ii$).
They are
$${\rm span} \{ Y_1,Y_2,Y_3 \}$$
where
$$Y_1 = X_1 - \alpha \tilde{X}^3, 
Y_2 = X_2 - \beta \tilde{X}^3, 
Y_3= X_3 + \alpha \tilde{X}^1 + \beta \tilde{X}^2,  \alpha, \beta \in \real,$$
 and the dual subalgebra w.r.t. $\langle \, . , . \rangle$ is in general
$${\rm span} \{ \tilde{Y}_1,\tilde{Y}_2,\tilde{Y}_3 \}$$
where 
$$\tilde{Y}_1 =  \gamma X_2 + \tilde{X}^1 - \gamma \beta \tilde{X}^3, \tilde{Y}_2 =
-\gamma X_1 + \tilde{X}^2 +\gamma \alpha \tilde{X}^3,  \tilde{Y}_3 = \tilde{X}^2, \, \gamma \in \real.$$
Structure coefficients in this new basis $Y_1,\ldots, \tilde{Y}_3  $ are identical with 
the original structure coefficients
for any 
$\alpha,\beta,\gamma$, therefore the Drinfeld double 
corresponding to ($7_0|2.ii$) allows no decomposition into other Manin triples and similarly for ($7_0|2.i$).

Concerning the Lie algebra structure, the Drinfeld doubles corresponding  to ($7_0|2.i$) and ($7_0|2.ii$) are isomorphic
as Lie algebras because they differ just by the sign of the bilinear form 
$\langle \ , \ \rangle$, and consequently the commutation relations implied by ad--invariance of 
$\langle \ , \ \rangle$ are the same. The other Drinfeld doubles specify different Lie algebras for the same reason as in the Section
\ref{nilp105}.

\subsection{Manin triples with the  Killing form of signature (0,0,6)}
\subsubsection{Case ${\rm \bf  dim} {\bf [\cd, \cd]=5}$}
This set contains Manin triples  ($7_{1}|1$),($7_{1}|2.i$) and ($7_{1}|2.ii$). 
They specify isomorphic Drinfeld doubles. For transformation matrices see subsection 
\ref{nilp105} and substitute $a=1$.

\subsubsection{Case ${ \rm \bf dim}{\bf [\cd, \cd]=3}$}
In this set, the only Manin triples that can lead to the same Drinfeld double are ($2|2.i$)
and ($2|2.ii$). To see that the Drinfeld doubles are different, it is sufficent to find the centres
$Z( \cd )$ of 
these Manin triples and restrict the form $\langle \,.,. \rangle$ to them. These restricted forms 
$\langle \,.,.\rangle_{Z( \cd )}$
have different signatures, therefore the Drinfeld doubles are non--isomorphic:
\begin{enumerate}
\item ($2|2.i$): $Z( \cd )={\rm span} \{ X_1, X_2 - \tilde{X}^2, \tilde{X}^3\}$, \\
signature of $\langle \,.,. \rangle_{Z( \cd )}=(0,1,2)$
\item ($2|2.ii$): $Z( \cd )={\rm span} \{ X_1, X_2 + \tilde{X}^2, \tilde{X}^3\}$, \\
signature of $\langle \,.,.\rangle_{Z( \cd )}=(1,0,2)$
\end{enumerate}
These Drinfeld doubles are isomorphic as Lie algebras because they differ just by the sign of the bilinear form 
$\langle \ , \ \rangle$ and the commutation relations are due to the ad--invariance the same.
\smallskip

\noindent {\bf Q.E.D.}

\section{Conclusions}
In this work we have constructed the complete list of six--dimensional real Drinfeld doubles up to their isomorphisms i.e. maps 
preserving both the Lie structure and an ad--invariant symmetric bilinear form $\langle , \rangle$ that define the double. 
The result is summarized in the theorem at the end of the Section \ref{sec3} and claims that there just 22 classes of 
the non--isomorphic Drinfeld doubles. Some of them contain one or two real parameters denoted $a$ and $b$. The number 22 is 
in a way conditional because e.g. the classes 9,14,18 could be united into one. 
The reason why they are given as separated 
classes is that they have different values of their invariants, namely the signature of the Killing form.

An important point that follows from the classification is that there are several different Drinfeld doubles corresponding to Lie algebras $so(1,3),\ sl(2,\real)\oplus sl(2,\real),\ sl(2,\real)\triangleright\real^3$ whereas on solvable Lie algebras the Drinfeld double is unique (in some cases up to the sign of the bilinear form). 
On the other hand there are Manin triples with one isotropic subalgebra abelian that are equivalent as Drinfeld doubles even though the other subalgebras are different (see $(6_0|1)$ and $(5|1)$ ). That's why it is necessary to investigate the (non)equivalence of the Manin triples of this form. Moreover the above given examples indicate the diversity of Drinfeld double structures one may encounter in higher dimensions.

Beside that from the present classification procedure one can find whether a given six--dimensional Lie algebra can be equipped by a suitable ad--invariant bilinear form and turned into a Drinfeld double (and how many such forms exist). The decisive aspects are the signature of the Killing form and the dimensions of the ideals $\cd_j,\cd^j$. The necessary condition is that they have the values occuring in the Table \ref{tab1}. The investigation then can be reduced to a direct check of equivalence with a particular six--dimensional Lie algebra (possibly after determination of abelian ideals and 
the factor algebras as in the Section \ref{nilp105}).

One can see that for many Drinfeld doubles there are several decompositions into Manin triples. For each Manin triple there is 
a pair of dual sigma models. Their equation of motion \cite{kli:pltd}
\be \langle\partial_\pm l l^{-1},{\cal E}^\pm\rangle=0 \lbl{eqm}\ee
are given by the Drinfeld double and a three--dimensional subspace ${\cal E}^+\subset\cd$ so that all these models (for fixed 
${\cal E}^+$) are equivalent. Moreover the scaling of $\langle , \rangle$ does not change the equations of motion \rf{eqm}) and consequently all the models corresponding to (non--isomorphic) Drinfeld doubles with different $b$ are equal. We can construct 
the explicit forms of the equations of motion for every Drinfeld double but without a physical motivation this does not make much sense.

Let us note that the complete sets of the equivalent sigma models for a fixed Drinfeld double are given by the so called modular 
space of the double. The construction of all non--isomorphic Manin triples for the double is the first step in the construction 
of the modular spaces.
\section{Appendix A: Bianchi algebras}
It is known that any 3--dimensional real Lie algebra can be brought to one of 11 forms by a change of basis.
These forms represent non--isomorphic Lie algebras and are conventionally known as Bianchi algebras.  
They are denoted by  ${\bf 1},\ldots,{\bf 5}$, ${\bf 6_a}$,${\bf 6_0}$,
${\bf 7_a}$,${\bf 7_0}$,${\bf 8}$,${\bf 9}$
 (see e.g. \cite {Landau}, in literature often uppercase 
roman numbers are used instead of arabic ones). The list of Bianchi algebras is given in decreasing order 
starting from simple algebras.

\begin{description}
\item[${\bf 9:}$] $[X_1,X_2]=X_3, \; [X_2,X_3] = X_1, \; [X_3,X_1] = X_2, $ (i.e. $so(3)$)
\item[${\bf 8:}$] $[X_1,X_2]=-X_3, \; [X_2,X_3] = X_1, \; [X_3,X_1] = X_2, $ (i.e. $sl(2,{\real})$)
\item[${\bf 7_a:}$] $[X_1,X_2]=-a X_2+X_3, \; [X_2,X_3] = 0, \; [X_3,X_1] = X_2+ a X_3, 
\; a >0, $
\item[${\bf 7_0:}$] $[X_1,X_2]=0, \; [X_2,X_3] = X_1, \; [X_3,X_1] = X_2, $
\item[${\bf 6_a:}$] $[X_1,X_2]=-a X_2-X_3, \; [X_2,X_3] = 0, \; [X_3,X_1] = X_2+ a X_3, 
\; a >0, \, a \neq 1 , $
\item[${\bf 6_0:}$] $[X_1,X_2]=0, \; [X_2,X_3] = X_1, \; [X_3,X_1] = -  X_2, $
\item[${\bf 5:}$] $[X_1,X_2]=-X_2, \; [X_2,X_3] = 0, \; [X_3,X_1] = X_3, $
\item[${\bf 4:}$] $[X_1,X_2]=-X_2+X_3, \; [X_2,X_3] = 0, \; [X_3,X_1] = X_3, $
\item[${\bf 3:}$] $[X_1,X_2]=-X_2-X_3, \; [X_2,X_3] = 0, \; [X_3,X_1] = X_2+X_3, $
\item[${\bf 2:}$] $[X_1,X_2]=0, \; [X_2,X_3] = X_1, \; [X_3,X_1] = 0,$
\item[${\bf 1:}$] $[X_1,X_2]=0, \; [X_2,X_3] = 0, \; [X_3,X_1] = 0 , $
\end{description}

One might use also another classification (used e.g. in \cite{gom:ctd}). In this notation the basis of the Lie algebra is usually 
written as $(e_0,e_1,e_2)$ and the classification is:
\begin{description}
\item[${\real}^3={\bf 1}$]: $[e_1,e_2]=0$, $[e_0,e_i]=0$
\item[$n_3={\bf 2}$]: $[e_1,e_2]=e_0$, $[e_0,e_i]=0$
\item[$r_3(\rho)$]: $[e_1,e_2]=0$, $[e_0,e_1]=e_1$, $[e_0,e_2]=\rho e_2$, $-1 \leq \rho \leq 1$ \\
This algebra is isomorphic to ${\bf 6_0}$ for $\rho =-1$, ${\bf 6_{\frac{\rho +1}{\rho-1}}}$ for $0 <|\rho| < 1$,
${\bf 3}$ for $\rho =0$ and ${\bf 5}$ for $\rho =1$.
\item[$r'_3(1)={\bf 4}$]: $[e_1,e_2]=0$, $[e_0,e_1]=e_1$, $[e_0,e_2]=e_1+ e_2$ 
\item[$s_3(\rho)$]: $[e_1,e_2]=0$, $[e_0,e_1]=\mu e_1-e_2$, $[e_0,e_2]=e_1+\mu e_2$, $\mu \geq 0$ \\
This algebra is isomorphic to ${\bf 7_0}$ for $\mu =0$ and ${\bf 7_\mu}$ for $\mu > 0 $.
\item[$sl(2,{\real})={\bf 8}$]
\item[$so(3)={\bf 9}$]
\end{description}
It is clear that this classification is more compact, on the other hand the classes in this classification contain 
algebras with different properties such as dimensions of commutant etc. and surprisingly the special cases of parameters 
we need to distinguish correspond in most cases to different Bianchi algebras. Therefore we use the Bianchi classification.

\section{Appendix B: 
List of Manin triples}
We present a list of Manin triples based on \cite{hlasno:mtriples}. 
The label of each Manin triple, e.g. 
${\bf  ( 8|5.ii|b)  }$, indicates the structure of the first subalgebra  $\cg$, 
e.g. Bianchi algebra ${\bf 8}$, the structure of the second subalgebra  $\tcg$, e.g. Bianchi algebra  
${\bf 5}$;
roman numbers $i$, $ii$ etc. (if present) distinguish between several possible 
pairings $\langle \, . , . \rangle$ of the subalgebras $\cg,\tcg$ and the parameter $b$ indicates
the Manin triples differing by the rescaling of $\langle \, . , . \rangle$
(if such Manin triples are not isomorphic).

The Lie structures of the subalgebras $\cg$ and $\tcg$ are written out in mutually dual bases 
$(X_1,X_2,X_3)$ and $(\tilde{X}^1,\tilde{X}^2,\tilde{X}^3)$ where the transformation (\ref{tfnb})
was used to bring $\cg$ to the standard Bianchi form. Because of (\ref{liebd})
this information specifies the Manin triple completely.

The dual Manin triples ($\cd,\tcg, \cg$) are not written explicitly
 but can be easily obtained by $X_j\leftrightarrow\tilde X^j$.

\begin{enumerate}
\item Manin triples with the first subalgebra $ \cg= {\bf 9}$:
\begin{description}
\item[${\bf (9|1 ) } $] : \\ 
  $[X_1,X_2]=X_3, \; [X_2,X_3] = X_1, \; [X_3,X_1] = X_2, $ \\
 $ [\tilde{X}^1,\tilde{X}^2]= 0 ,  \, 
[\tilde{X}^2,\tilde{X}^3] = 0 , \; 
[\tilde{X}^3,\tilde{X}^1] = 0. $ 
\vskip4mm  \item[${\bf (9|5|b )}$]
  $[X_1,X_2]=X_3, \; [X_2,X_3] = X_1, \; [X_3,X_1] = X_2, $ \\
  $ \; \; [\tilde{X}^1,\tilde{X}^2]=- b \tilde{X}^2, \, 
[\tilde{X}^2,\tilde{X}^3] = 0 , \; 
[\tilde{X}^3,\tilde{X}^1] = b \tilde{X}^3, \; b > 0  . $ 
\end{description}

\vskip4mm  \item Manin triples with the first subalgebra $ \cg= {\bf 8}$:
\begin{description} 
\vskip4mm  \item[${\bf ( 8|1) } $] : \\
$  [X_1,X_2]=-X_3, \; [X_2,X_3] = X_1, \; [X_3,X_1] = X_2, $
\\v
$[\tilde{X}^1,\tilde{X}^2]= 0 ,  \, 
[\tilde{X}^2,\tilde{X}^3] = 0 , \; 
[\tilde{X}^3,\tilde{X}^1] = 0. $ 
\vskip4mm  \item[${\bf ( 8|5.i|b) } $] : \\
$  [X_1,X_2]=-X_3, \; [X_2,X_3] = X_1, \; [X_3,X_1] = X_2, $
\\
$ [\tilde{X}^1,\tilde{X}^2]=- b \tilde{X}^2, \, 
[\tilde{X}^2,\tilde{X}^3] = 0 , \; 
[\tilde{X}^3,\tilde{X}^1] = b \tilde{X}^3, \; b >0 . $ 
\vskip4mm  \item[${\bf ( 8|5.ii|b) } $] : \\
$  [X_1,X_2]=-X_3, \; [X_2,X_3] = X_1, \; [X_3,X_1] = X_2, $
\\
$ [\tilde{X}^1,\tilde{X}^2]=0, \, 
[\tilde{X}^2,\tilde{X}^3] = b  \tilde{X}^2 , \; 
[\tilde{X}^3,\tilde{X}^1] = - b \tilde{X}^1, \; b >0 . $ 
\vskip4mm  \item[${\bf ( 8|5.iii) } $] : \\
$  [X_1,X_2]=-X_3, \; [X_2,X_3] = X_1, \; [X_3,X_1] = X_2, $
\\
$ [\tilde{X}^1,\tilde{X}^2]= \tilde{X}^2, \, 
[\tilde{X}^2,\tilde{X}^3] =   \tilde{X}^2  , \; 
[\tilde{X}^3,\tilde{X}^1] = - ( \tilde{X}^1+\tilde{X}^3 ) . $
\end{description}

\vskip4mm  \item Manin triples with the first subalgebra $ \cg= {\bf 7_{a}}$:
\begin{description} 
\vskip4mm  \item[${\bf ( 7_a|1) } $] : \\
$  [X_1,X_2]=-a X_2+X_3, \; [X_2,X_3] = 0, \; [X_3,X_1] = X_2+ a X_3, 
\; a >0, $
\\
$  [\tilde{X}^1,\tilde{X}^2]= 0, \, 
[\tilde{X}^2,\tilde{X}^3] = 0 , \; 
[\tilde{X}^3,\tilde{X}^1] = 0. $ 
\vskip4mm  \item[${\bf ( 7_a|2.i) } $] : \\
$  [X_1,X_2]=-a X_2+X_3, \; [X_2,X_3] = 0, \; [X_3,X_1] = X_2+ a X_3, 
\; a >0, $
\\
$ [\tilde{X}^1,\tilde{X}^2]=0, \, 
[\tilde{X}^2,\tilde{X}^3] = \tilde{X}^1 , \; 
[\tilde{X}^3,\tilde{X}^1] = 0 .$ 
\vskip4mm  \item[${\bf ( 7_a|2.ii) } $] : \\
$  [X_1,X_2]=-a X_2+X_3, \; [X_2,X_3] = 0, \; [X_3,X_1] = X_2+ a X_3, 
\; a >0, $
\\
$ [\tilde{X}^1,\tilde{X}^2]=0, \, 
[\tilde{X}^2,\tilde{X}^3] = - \tilde{X}^1 , \; 
[\tilde{X}^3,\tilde{X}^1] = 0 . $
\vskip4mm  \item[${\bf ( 7_a|7_{1/a}|b) } $] : \\ 
$ [X_1,X_2]=-a X_2+X_3, \; [X_2,X_3] = 0, \; [X_3,X_1] = X_2+ a X_3, 
\; a >0, $
\\
$ [\tilde{X}^1,\tilde{X}^2]= b ( -  \frac{1}{a} \tilde{X}^2 + \tilde{X}^3), \, 
[\tilde{X}^2,\tilde{X}^3] = 0, \; 
[\tilde{X}^3,\tilde{X}^1] = b (\tilde{X}^2+ \frac{1}{a} \tilde{X}^3), 
\; b \in { \real} - \{ 0 \} . $ 
\end{description}

\vskip4mm  \item Manin triples with the first subalgebra $ \cg= {\bf 7_0}$:
\begin{description} 
\vskip4mm  \item[${\bf ( 7_0|1) } $] : \\
$  [X_1,X_2]=0, \; [X_2,X_3] = X_1, \; [X_3,X_1] = X_2, 
 $
\\
$ [\tilde{X}^1,\tilde{X}^2]= 0, \, 
[\tilde{X}^2,\tilde{X}^3] = 0 , \; 
[\tilde{X}^3,\tilde{X}^1] = 0. $ 
\vskip4mm  \item[${\bf ( 7_0|2.i) } $] : \\
$  [X_1,X_2]=0, \; [X_2,X_3] = X_1, \; [X_3,X_1] = X_2, $
\\
$ [\tilde{X}^1,\tilde{X}^2]= \tilde{X}^3, \, 
[\tilde{X}^2,\tilde{X}^3] = 0 , \; 
[\tilde{X}^3,\tilde{X}^1] = 0 .$ 
\vskip4mm  \item[${\bf ( 7_0|2.ii) } $] : \\
$  [X_1,X_2]=0, \; [X_2,X_3] = X_1, \; [X_3,X_1] = X_2, $
\\
$ [\tilde{X}^1,\tilde{X}^2]= - \tilde{X}^3, \, 
[\tilde{X}^2,\tilde{X}^3] = 0 , \; 
[\tilde{X}^3,\tilde{X}^1] = 0 . $
\vskip4mm  \item[${\bf ( 7_0|4|b) } $] : \\
$  [X_1,X_2]=0, \; [X_2,X_3] = X_1, \; [X_3,X_1] = X_2, $
\\
$ [\tilde{X}^1,\tilde{X}^2]= b ( - \tilde{X}^2 + \tilde{X}^3), \, 
[\tilde{X}^2,\tilde{X}^3] = 0, \; 
[\tilde{X}^3,\tilde{X}^1] = b  \tilde{X}^3, 
\; b \in { \real} - \{ 0 \} . $ 
\vskip4mm  \item[${\bf ( 7_0|5.i) } $] : \\
$  [X_1,X_2]=0, \; [X_2,X_3] = X_1, \; [X_3,X_1] = X_2, $
\\
$ [\tilde{X}^1,\tilde{X}^2]=  - \tilde{X}^2 , \, 
[\tilde{X}^2,\tilde{X}^3] = 0, \; 
[\tilde{X}^3,\tilde{X}^1] =  \tilde{X}^3, 
. $ 
\vskip4mm  \item[${\bf ( 7_0|5.ii|b) } $] : \\
$  [X_1,X_2]=0, \; [X_2,X_3] = X_1, \; [X_3,X_1] = X_2, $
\\
$  [\tilde{X}^1,\tilde{X}^2]=  0 , \, 
[\tilde{X}^2,\tilde{X}^3] = b \tilde{X}^2, \; 
[\tilde{X}^3,\tilde{X}^1] = -b \tilde{X}^1, 
\, b >0  . $
\end{description}

\vskip4mm  \item Manin triples with the first subalgebra $ \cg= {\bf 6_{a}}$:
\begin{description} 
\vskip4mm  \item[${\bf ( 6_{a}|1) } $] : \\
$  [X_1,X_2]=-a X_2-X_3, \; [X_2,X_3] = 0, \; [X_3,X_1] = X_2+ a X_3, 
\; a >0, \, a \neq 1 , $

$  [\tilde{X}^1,\tilde{X}^2]= 0, \, 
[\tilde{X}^2,\tilde{X}^3] = 0 , \; 
[\tilde{X}^3,\tilde{X}^1] = 0. $ 
\vskip4mm  \item[${\bf ( 6_{a}|2) } $] : \\ 
$  [X_1,X_2]=-a X_2-X_3, \; [X_2,X_3] = 0, \; [X_3,X_1] = X_2+ a X_3, 
\; a >0, \, a \neq 1 , $

$ [\tilde{X}^1,\tilde{X}^2]=0, \, 
[\tilde{X}^2,\tilde{X}^3] = \tilde{X}^1 , \; 
[\tilde{X}^3,\tilde{X}^1] = 0 . $ 
\vskip4mm  \item[${\bf ( 6_{a}|6_{1/a}.i|b) } $] : \\
$  [X_1,X_2]=-a X_2-X_3, \; [X_2,X_3] = 0, \; [X_3,X_1] = X_2+ a X_3, 
\; a >0, \, a \neq 1 , $
\\
$ [\tilde{X}^1,\tilde{X}^2]=- b ( \frac{1}{a} \tilde{X}^2+\tilde{X}^3), \, 
[\tilde{X}^2,\tilde{X}^3] = 0, \; 
[\tilde{X}^3,\tilde{X}^1] = b (\tilde{X}^2+ \frac{1}{a} \tilde{X}^3), \; b \in { \real} - \{ 0 \} . $ 
\vskip4mm  \item[${\bf ( 6_{a}|6_{1/a}.ii) } $] : \\
$  [X_1,X_2]=-a X_2-X_3, \; [X_2,X_3] = 0, \; [X_3,X_1] = X_2+ a X_3, 
\; a >0, \, a \neq 1 , $
\\
$ [\tilde{X}^1,\tilde{X}^2]= \tilde{X}^1 , \, 
[\tilde{X}^2,\tilde{X}^3] = \frac{a+1}{a-1} (\tilde{X}^2+\tilde{X}^3) , \; 
[\tilde{X}^3,\tilde{X}^1] = \tilde{X}^1 . $
\vskip4mm  \item[${\bf ( 6_{a}|6_{1/a}.iii) } $] : \\
$  [X_1,X_2]=-a X_2-X_3, \; [X_2,X_3] = 0, \; [X_3,X_1] = X_2+ a X_3, 
\; a >0, \, a \neq 1 , $
\\
$ [\tilde{X}^1,\tilde{X}^2]= \tilde{X}^1 , \, 
[\tilde{X}^2,\tilde{X}^3] = \frac{a-1}{a+1} (-\tilde{X}^2+\tilde{X}^3) , \; 
[\tilde{X}^3,\tilde{X}^1] = -\tilde{X}^1 .$ 
\end{description}

\vskip4mm  \item Manin triples with the first subalgebra $ \cg= {\bf 6_0}$:
\begin{description}
\vskip4mm  \item[${\bf ( 6_0|1) } $] : \\ 
$  [X_1,X_2]=0, \; [X_2,X_3] = X_1, \; [X_3,X_1] = -  X_2, $
\\
$ [\tilde{X}^1,\tilde{X}^2]= 0, \, 
[\tilde{X}^2,\tilde{X}^3] = 0 , \; 
[\tilde{X}^3,\tilde{X}^1] = 0. $ 
\vskip4mm  \item[${\bf ( 6_0|2) } $] : \\ 
$  [X_1,X_2]=0, \; [X_2,X_3] = X_1, \; [X_3,X_1] = -  X_2, $
\\
 $  [\tilde{X}^1,\tilde{X}^2]= \tilde{X}^3, \, 
[\tilde{X}^2,\tilde{X}^3] = 0 , \; 
[\tilde{X}^3,\tilde{X}^1] = 0 .$ 
\vskip4mm  \item[${\bf ( 6_0|4.i|b) } $] : \\ 
$  [X_1,X_2]=0, \; [X_2,X_3] = X_1, \; [X_3,X_1] = -  X_2, $
\\
$ [\tilde{X}^1,\tilde{X}^2]= b ( - \tilde{X}^2 + \tilde{X}^3), \, 
[\tilde{X}^2,\tilde{X}^3] = 0, \; 
[\tilde{X}^3,\tilde{X}^1] = b  \tilde{X}^3, 
\; b \in { \real} - \{ 0 \} . $
\vskip4mm  \item[${\bf ( 6_0|4.ii) } $] : \\ 
$  [X_1,X_2]=0, \; [X_2,X_3] = X_1, \; [X_3,X_1] = -  X_2, $
\\
$ [\tilde{X}^1,\tilde{X}^2]=  ( - \tilde{X}^1 + \tilde{X}^2 + \tilde{X}^3), \, 
[\tilde{X}^2,\tilde{X}^3] = \tilde{X}^3, \; 
[\tilde{X}^3,\tilde{X}^1] =  - \tilde{X}^3 
 . $
\vskip4mm  \item[${\bf ( 6_0|5.i) } $] : \\  
$  [X_1,X_2]=0, \; [X_2,X_3] = X_1, \; [X_3,X_1] = -  X_2, $
\\
$ [\tilde{X}^1,\tilde{X}^2]=  - \tilde{X}^2 , \, 
[\tilde{X}^2,\tilde{X}^3] = 0, \; 
[\tilde{X}^3,\tilde{X}^1] =  \tilde{X}^3. $ 
\vskip4mm  \item[${\bf ( 6_0|5.ii) } $] : \\  
$  [X_1,X_2]=0, \; [X_2,X_3] = X_1, \; [X_3,X_1] = -  X_2, $
\\
$ [\tilde{X}^1,\tilde{X}^2]=  - \tilde{X}^1+ \tilde{X}^2 , \, 
[\tilde{X}^2,\tilde{X}^3] =  \tilde{X}^3, \; 
[\tilde{X}^3,\tilde{X}^1] = - \tilde{X}^3. $ 
\vskip4mm  \item[${\bf ( 6_0|5.iii|b) } $] : \\  
$  [X_1,X_2]=0, \; [X_2,X_3] = X_1, \; [X_3,X_1] = -  X_2, $
\\
$ [\tilde{X}^1,\tilde{X}^2]=  0 , \, 
[\tilde{X}^2,\tilde{X}^3] = - b\tilde{X}^2, \; 
[\tilde{X}^3,\tilde{X}^1] =  b\tilde{X}^1, \ b>0 . $ 
\end{description}

\vskip4mm  \item Manin triples with the first subalgebra $ \cg= {\bf 5}$:
\begin{description} 
\vskip4mm  \item[${\bf ( 5|1) } $] : \\ 
$  [X_1,X_2]=-X_2, \; [X_2,X_3] = 0, \; [X_3,X_1] = X_3, $
\\
$ [\tilde{X}^1,\tilde{X}^2]= 0, \, 
[\tilde{X}^2,\tilde{X}^3] = 0 , \; 
[\tilde{X}^3,\tilde{X}^1] = 0. $ 
\vskip4mm  \item[${\bf ( 5|2.i) } $] : \\   
$  [X_1,X_2]=-X_2, \; [X_2,X_3] = 0, \; [X_3,X_1] = X_3, $
\\
$  [\tilde{X}^1,\tilde{X}^2]= 0, \, 
[\tilde{X}^2,\tilde{X}^3] = \tilde{X}^1 , \; 
[\tilde{X}^3,\tilde{X}^1] =  0  .$
\vskip4mm  \item[${\bf ( 5|2.ii) } $] : \\
$  [X_1,X_2]=-X_2, \; [X_2,X_3] = 0, \; [X_3,X_1] = X_3, $
\\
$ [\tilde{X}^1,\tilde{X}^2]= \tilde{X}^3, \, 
[\tilde{X}^2,\tilde{X}^3] = 0 , \; 
[\tilde{X}^3,\tilde{X}^1] =  0  .$
\end{description}
  and dual Manin triples ${\bf ( \cg \leftrightarrow \tcg ) } $ to Manin triples given above for  
$\cg = {\bf 6_0}$, ${\bf 7_0} $, ${\bf 8}$, ${\bf 9}$.   

\vskip4mm  \item Manin triples with the first subalgebra $ \cg= {\bf 4}$:
\begin{description} 
\vskip4mm  \item[${\bf ( 4|1) } $] : \\ 
$  [X_1,X_2]=-X_2+X_3, \; [X_2,X_3] = 0, \; [X_3,X_1] = X_3, $
\\
$ [\tilde{X}^1,\tilde{X}^2]= 0, \, 
[\tilde{X}^2,\tilde{X}^3] = 0 , \; 
[\tilde{X}^3,\tilde{X}^1] = 0. $ 
\vskip4mm  \item[${\bf ( 4|2.i) } $] : \\ 
$  [X_1,X_2]=-X_2+X_3, \; [X_2,X_3] = 0, \; [X_3,X_1] = X_3, $
\\
 $ [\tilde{X}^1,\tilde{X}^2]= 0, \, 
[\tilde{X}^2,\tilde{X}^3] = \tilde{X}^1 , \; 
[\tilde{X}^3,\tilde{X}^1] =  0  .$
\vskip4mm  \item[${\bf ( 4|2.ii) } $] : \\
$  [X_1,X_2]=-X_2+X_3, \; [X_2,X_3] = 0, \; [X_3,X_1] = X_3, $
\\
 $ [\tilde{X}^1,\tilde{X}^2]= 0, \, 
[\tilde{X}^2,\tilde{X}^3] = - \tilde{X}^1 , \; 
[\tilde{X}^3,\tilde{X}^1] =  0  .$
\vskip4mm  \item[${\bf ( 4|2.iii|b) } $] : \\ 
$  [X_1,X_2]=-X_2+X_3, \; [X_2,X_3] = 0, \; [X_3,X_1] = X_3, $
\\
 $ [\tilde{X}^1,\tilde{X}^2]= 0, \, 
[\tilde{X}^2,\tilde{X}^3] = 0 , \; 
[\tilde{X}^3,\tilde{X}^1] = b \tilde{X}^2, \, b \in \real- \{ 0 \} .$ 
\end{description}
  and dual Manin triples ${\bf ( \cg \leftrightarrow \tcg ) } $ to Manin triples given above for 
$\cg = {\bf 6_0}$, ${\bf 7_0}$.  

\vskip4mm  \item Manin triples with the first subalgebra $ \cg= {\bf 3}$:
\begin{description} 
\vskip4mm  \item[${\bf ( 3|1) } $] : \\
$  [X_1,X_2]=-X_2-X_3, \; [X_2,X_3] = 0, \; [X_3,X_1] = X_2+X_3, $
\\
$ [\tilde{X}^1,\tilde{X}^2]= 0, \, 
[\tilde{X}^2,\tilde{X}^3] = 0 , \; 
[\tilde{X}^3,\tilde{X}^1] = 0. $ 
\vskip4mm  \item[${\bf ( 3|2) } $] : \\ 
$  [X_1,X_2]=-X_2-X_3, \; [X_2,X_3] = 0, \; [X_3,X_1] = X_2+X_3, $
\\
$ [\tilde{X}^1,\tilde{X}^2]=0, \, 
[\tilde{X}^2,\tilde{X}^3] = \tilde{X}^1 , \; 
[\tilde{X}^3,\tilde{X}^1] = 0 . $ 
\vskip4mm  \item[${\bf ( 3|3.i) } $] : \\
$  [X_1,X_2]=-X_2-X_3, \; [X_2,X_3] = 0, \; [X_3,X_1] = X_2+X_3, $
\\
$ [\tilde{X}^1,\tilde{X}^2]=- b (\tilde{X}^2+\tilde{X}^3), \, 
[\tilde{X}^2,\tilde{X}^3] = 0, \; 
[\tilde{X}^3,\tilde{X}^1] = b (\tilde{X}^2+\tilde{X}^3), \; b \in { \real} - \{ 0 \} . $ 
\vskip4mm  \item[${\bf ( 3|3.ii) } $] : \\
$  [X_1,X_2]=-X_2-X_3, \; [X_2,X_3] = 0, \; [X_3,X_1] = X_2+X_3, $
\\
$ [\tilde{X}^1,\tilde{X}^2]= 0 , \, 
[\tilde{X}^2,\tilde{X}^3] = \tilde{X}^2+\tilde{X}^3 , \; 
[\tilde{X}^3,\tilde{X}^1] =  0 . $
\vskip4mm  \item[${\bf ( 3|3.iii) } $] : \\
$  [X_1,X_2]=-X_2-X_3, \; [X_2,X_3] = 0, \; [X_3,X_1] = X_2+X_3, $
\\
$ [\tilde{X}^1,\tilde{X}^2]= \tilde{X}^1 , \, 
[\tilde{X}^2,\tilde{X}^3] = 0 , \; 
[\tilde{X}^3,\tilde{X}^1] = - \tilde{X}^1 . $
\end{description}

\vskip4mm  \item Manin triples with the first subalgebra $ \cg= {\bf 2}$:
\begin{description} 
\vskip4mm  \item[${\bf ( 2|1) } $] : \\ 
$  [X_1,X_2]=0, \; [X_2,X_3] = X_1, \; [X_3,X_1] = 0,$
\\
$ [\tilde{X}^1,\tilde{X}^2]= 0, \, 
[\tilde{X}^2,\tilde{X}^3] = 0 , \; 
[\tilde{X}^3,\tilde{X}^1] = 0. $ 
\vskip4mm  \item[${\bf ( 2|2.i) } $] : \\
$  [X_1,X_2]=0, \; [X_2,X_3] = X_1, \; [X_3,X_1] = 0,
$
\\
 $ [\tilde{X}^1,\tilde{X}^2]= \tilde{X}^3, \, 
[\tilde{X}^2,\tilde{X}^3] = 0 , \; 
[\tilde{X}^3,\tilde{X}^1] =  0  .$
\vskip4mm  \item[${\bf ( 2|2.ii) } $] : \\
$  [X_1,X_2]=0, \; [X_2,X_3] = X_1, \; [X_3,X_1] = 0,
$
\\
 $ [\tilde{X}^1,\tilde{X}^2]= -\tilde{X}^3, \, 
[\tilde{X}^2,\tilde{X}^3] = 0 , \; 
[\tilde{X}^3,\tilde{X}^1] =  0  .$
\end{description}
 and dual Manin triples ${\bf ( \cg \leftrightarrow \tcg ) } $ to Manin triples given above  for 
$\cg= {\bf 3}$, ${\bf 4}$, ${\bf 6_0}$, ${\bf 6_a}$,  ${\bf 7_0}$, ${\bf 7_a}$.    

\vskip4mm  \item Manin triples with the first subalgebra $ \cg= {\bf 1}$:
\begin{description}
\vskip4mm 
 \item[${\bf ( 1|1) } $] : \\
$  [X_1,X_2]=0, \; [X_2,X_3] = 0, \; [X_3,X_1] = 0 , $
\\
$ [\tilde{X}^1,\tilde{X}^2]= 0, \, 
[\tilde{X}^2,\tilde{X}^3] = 0 , \; 
[\tilde{X}^3,\tilde{X}^1] = 0. $ 
\end{description}
and dual Manin triples ${\bf ( \cg \leftrightarrow \tcg ) } $ to Manin triples given above  for 
$\cg= {\bf 2}$--${\bf 9}$.

\end{enumerate}

\end{document}